\documentclass[notitlepage,11pt]{article}
\usepackage{amssymb,amsmath
}

\catcode`\@=11
\@addtoreset{equation}{section}

\catcode`\@=12

\usepackage{color}

\def\R{\mathbb{R}}

\def\f{\varphi}
\def\irn{\int\limits_{\R^n}}

\def\eps{\varepsilon}
\def\div{{\rm div}}
\def\mstar{{2_m^*}}

\def\tildeSobolev{\widetilde{H}^m(\Omega)}

\def\SoRn{H^m(\R^n)}
\def\SoRn{H^m(\mathbb R^n)}

\def\proof{\noindent{\textbf{Proof. }}}
\def\QED{\hfill {$\square$}\goodbreak \medskip}

\newtheorem{Theorem}{Theorem}
{Lemma}
\newtheorem{Proposition}
{Proposition}
{Corollary}
{Remark}
{Definition}

\begin{document}

\title 
{On fractional Laplacians -- 3}

\author{Roberta Musina\footnote{Dipartimento di Matematica ed Informatica, Universit\`a di Udine,
via delle Scienze, 206 -- 33100 Udine, Italy. Email: {roberta.musina@uniud.it}. 
{Partially supported by Miur-PRIN 2009WRJ3W7-001 ``Fenomeni di concentrazione e {pro\-ble\-mi} di analisi geometrica''.}}~ and
Alexander I. Nazarov\footnote{
St.Petersburg Department of Steklov Institute, Fontanka 27, St.Petersburg, 191023, Russia, 
and St.Petersburg State University, 
Universitetskii pr. 28, St.Petersburg, 198504, Russia. E-mail: al.il.nazarov@gmail.com.
Supported by RFBR grant 14-01-00534 and by
St.Petersburg University grant 6.38.670.2013.}
}

\date{}

\maketitle

\begin{abstract}
We investigate the role of the noncompact group of dilations in $\R^n$ on the difference
of the quadratic forms associated to the fractional Dirichlet and Navier Laplacians.
Then we apply our results to study the Brezis--Nirenberg effect in
two families of noncompact boundary value problems involving the Navier-Laplacian . 
\footnotesize

\medskip

\medskip

\noindent
\textbf{2010 Mathematics Subject Classfication:} 47A63; 35A23.
\end{abstract}

\normalsize

\bigskip

\section{Introduction}
\label{S:Introduction}

The Sobolev space $\SoRn=W^m_2(\mathbb{R}^n)$, $m\in\mathbb R$, is the space of distributions $u\in{\cal S}'(\mathbb R^n)$
with finite norm
$$
\|u\|_m^2=\int\limits_{\mathbb R^n}\left(1+|\xi|^2\right)^m|\mathcal Fu(\xi)|^2\,d\xi,
$$
see for instance Section~2.3.3 of the monograph \cite{Tr}. Here ${\cal F}$ denotes
the Fourier transform
$$
\mathcal F{u}(\xi)=\frac{1}{(2\pi)^{n/2}}\int\limits_{\mathbb R^n} e^{-i\,\!\xi\cdot x}u(x)\,dx.
$$

For arbitrary $m\in\mathbb R$ we define fractional Laplacian on $\mathbb R^n$ by the quadratic 
form 
$$
Q_m[u]=((-\Delta)^mu,u):=\int\limits_{\mathbb R^n}|\xi|^{2m}|{\cal F}u(\xi)|^2 d\xi,
$$
with domain
$$
{\rm Dom}(Q_m)=
\{u\in {\cal S}'(\mathbb{R}^n) \,:\, Q_m[u]<\infty\}.
$$

Let $\Omega$ be a bounded and smooth domain in $\R^n$. 
%
We introduce the ``Dirichlet'' fractional Laplacian in $\Omega$ (denoted by $(-\Delta_{\Omega})^m_D$) as 
the restriction of $(-\Delta)^m$. The domain of its quadratic form is
$$
{\rm Dom}(Q_{m,\Omega}^D)=\{u\in {\rm Dom}(Q_m) \,:\, {\rm supp}\, u\subset\overline{\Omega}\}.
$$
Also we define the ``Navier'' fractional Laplacian as the $m$-th power of the conventional
Dirichlet Laplacian in the sense of spectral theory. Its quadratic form reads
$$
Q_{m,\Omega}^N[u]=((-\Delta_{\Omega})^m_Nu,u):=\sum\nolimits_j\lambda_j^m\cdot|(u,\varphi_j)|^2.
$$
Here, $\lambda_j$ and $\varphi_j$ are eigenvalues and eigenfunctions of the Dirichlet Laplacian in $\Omega$, respectively,
and ${\rm Dom}(Q_{m,\Omega}^N)$ consists of distributions in $\Omega$ such that $Q_{m,\Omega}^N[u]<\infty$.

For $m=1$ these operators evidently coincide: $(-\Delta_{\Omega})_N=(-\Delta_{\Omega})_D$.
We emphasize that, in contrast to $(-\Delta_{\Omega})^m_N$, the operator $(-\Delta_{\Omega})^m_D$ is not the $m$-th power of the 
Dirichlet Laplacian for $m\ne1$. 

It is well known that for $m>0$ quadratic forms $Q_{m,\Omega}^D$ and $Q_{m,\Omega}^N$ generate Hilbert 
structures on their domains, and 
$$
{\rm Dom}(Q_{m,\Omega}^D)=\tildeSobolev\subseteq{\rm Dom}(Q_{m,\Omega}^N),
$$
where
$$\tildeSobolev=\{u\in \SoRn\,:\,{\rm supp}\, u\subset\overline{\Omega}\}~\!.
$$
It is also easy to see that for $m\in\mathbb{N}$, $u\in\tildeSobolev$
$$
Q_{m,\Omega}^D[u]=Q_{m,\Omega}^N[u].
$$

In \cite{FL} and \cite{FL2} we compared the operators $(-\Delta_{\Omega})^m_D$ and $(-\Delta_{\Omega})^m_N$ 
for non-integer $m$. It turned out that the difference between their quadratic forms is positive or negative 
depending on the fact whether $\lfloor m\rfloor$ is odd or even.  However, roughly speaking, this difference 
disappears as $\Omega\to \R^n$.

Namely, denote by $F(\Omega)$ the class of smooth and bounded domains containing $\Omega$.
For any $u\in {\rm Dom}(Q_{m,\Omega}^D)$ the form $Q_{m,\Omega'}^D[u]$ does not depend on 
$\Omega'\in F(\Omega)$ while the form $Q_{m,\Omega'}^N[u]$ does depend on $\Omega'\supset\Omega$, 
and the following relations hold.

\begin{Proposition}
\label{P:old} (\cite[Theorem 2]{FL2}).
Let $m>-1$, $m\notin\mathbb{N}_0$. If $u\in {\rm Dom}(Q_{m,\Omega}^D)$, then
\begin{eqnarray}
Q_{m,\Omega}^D[u]=\inf_{\Omega'\in F(\Omega)}Q_{m,\Omega'}^N[u], & {\rm if} & 2k<m<2k+1,\ \ k\in\mathbb N_0;
\label{old1}\\
Q_{m,\Omega}^D[u]=\sup_{\Omega'\in F(\Omega)}Q_{m,\Omega'}^N[u], & {\rm if} & 2k-1<m<2k,\ \ k\in\mathbb N_0.
\label{old2}
\end{eqnarray}
\end{Proposition}

The main result of our paper is a quantitative version of Proposition \ref{P:old}. 

\begin{Theorem}
\label{T:main0}
Assume that $m>0$, $m\notin\mathbb{N}$. Let $u\in \tildeSobolev$, and let 
${\rm supp}(u)\subset B_r\subset B_R\subset\Omega$. Then
\begin{eqnarray}
Q_{m,\Omega}^N[u]\le Q_{m,\Omega}^D[u]+\frac {{C(n,m)}~\! R^n}{(R-r)^{2n+2m}}\cdot\|u\|^2_{L_1(\Omega)}, 
& {\rm if} & \lfloor m\rfloor\vdots\, 2;
\label{new1}\\
Q_{m,\Omega}^D[u]\le Q_{m,\Omega}^N[u]+\frac {{C(n,m)}~\! R^n}{(R-r)^{2n+2m}}\cdot\|u\|^2_{L_1(\Omega)}, 
& {\rm if} & \lfloor m\rfloor\!{\not\vdots}\, 2.
\label{new2}
\end{eqnarray}
\end{Theorem}

The proof of Theorem \ref{T:main0} is given in Section \ref{proof}. In Section \ref{BNN}
we apply this result for studying the equations\footnote{we assume that $0\in\Omega$.}
\begin{gather}
\label{eq:problem1}
(-\Delta_{\Omega})^m_Nu=\lambda(-\Delta_{\Omega})^s_N u+|u|^{\mstar-2}u\quad\textrm{in $\Omega$,}\\
\nonumber
~\\
\label{eq:problem2}
(-\Delta_{\Omega})^m_Nu=\lambda|x|^{-2s} u+|u|^{\mstar-2}u\quad\textrm{in $\Omega$,}
\end{gather}
where $0\le s<m<\frac{n}{2}$ and $\mstar=\frac{2n}{n-2m}$. By solution of (\ref{eq:problem1}) or 
(\ref{eq:problem2}) we mean a weak solution from ${\rm Dom}(Q_{m,\Omega}^N)$.

In the basic paper \cite{BN} by Brezis and Nirenberg a remarkable phenomenon was discovered
for the problem
\begin{equation}
\label{eq:BN_problem}
-\Delta u=\lambda u+|u|^{\frac{4}{n-2}}u\quad\textrm{in $\Omega$,}\qquad
u=0\quad
\textrm{on $\partial\Omega$,}
\end{equation}
which coincides with (\ref{eq:problem1}) and (\ref{eq:problem2}) with $n>2$, $m=1$, $s=0$.
Namely, the existence of a nontrivial solution for any small $\lambda>0$ holds if $n\ge 4$; in contrast, 
for $n=3$ non-existence phenomena for any sufficiently small $\lambda>0$ can be observed.
For this reason, the dimension $n=3$ has been  named {\em critical}
for problem (\ref{eq:BN_problem}) (compare with \cite{PuSe}, \cite{GGS}).

As was pointed out in \cite{MN-BN}, the Brezis--Nirenberg effect is a nonlinear analog of the so-called 
zero-energy resonance for the Schr\"odinger operators (see, e.g., \cite{Ya1} and
\cite[pp.287--288]{Ya}).

After \cite{BN}, a large number of papers have been focussed
on studying the effect of lower order linear perturbations in noncompact variational problems,
see for instance the list of references included in \cite[Chapter 7]{GGS} about the case
$m\in \mathbb N$, $s=0$, and the recent paper 
\cite{MN-BN}, where a survey of earlier results for the Dirichlet case
was given. For the Navier case with non-integer $m$, the only papers we know consider 
$m\in(0,1)$ and $s=0$, see 
\cite{Tan} and \cite{BaCPS}. See also the recent paper \cite{CdPS} and references therein 
for nonlinear lower-order perturbations.

We consider {the} general case and prove the following result (see Section \ref{BNN} for a more precise 
statement), {that corresponds to} 
\cite[Theorem 4.2]{MN-BN} .

\begin{Theorem}\label{main1}
Let $0\le s<m<\frac n2$. If $s\ge 2m-\frac{n}{2}$ then $n$ is not a critical dimension for the 
(\ref{eq:problem1}) and (\ref{eq:problem2}). This means that both these equations have ground 
state solutions for all sufficiently small $\lambda>0$.
\end{Theorem}

Let us recall some notation. $B_R$ is the ball with radius $R$ centered at the origin, $\mathbb S_R$ is
its boundary
.
We denote by $c$ with indices all explicit constants while $C$ without indices stand
for all inessential positive constants. To indicate that $C$ depends on some parameter $a$, we
write $C(a)$.

\section{Proof of Theorem \ref{T:main0}}\label{proof}

Notice that we can assume $u\in {\cal C}^\infty_0(\Omega)$, the general case is covered by approximation.
\smallskip

\noindent 
{\bf Proof of (\ref{new1})}. Let $m=2k+\sigma$, $k\in\mathbb N_0$, $\sigma\in(0,1)$. Denote by $w^D(x,y)$, $x\in\mathbb{R}^n$, 
$y>0$, the Caffarelli--Silvestre extension of $(-\Delta)^ku$ (see \cite{CaSi}), that is the solution of the 
boundary value problem
$$-\div (y^{1-2\sigma }\nabla w)=0\quad \mbox{in}\quad \mathbb R^n\times\mathbb R_+;
\qquad w\big|_{y=0}=(-\Delta)^ku,
$$
given by the generalized Poisson formula
\begin{equation}\label{Poisson}
w^D(x,y)=c_1(n,\sigma) \irn\frac{y^{2\sigma}\,(-\Delta)^k u(\xi)}
{\left(|x-\xi|^2+{y^2}\right)^{\frac{n+2\sigma}{2}}}\,d\xi .
\end{equation}
{In \cite{CaSi} it is also proved that}
\begin{equation}\label{QD}
Q_{m,\Omega}^D[u]=Q_{\sigma,\Omega}^D[(-\Delta)^ku]=c_2(n,\sigma)\int\limits_0^\infty\int\limits_{\mathbb R^n} 
y^{1-2\sigma}|\nabla w^D|^2\,dxdy.
\end{equation}

Integrating by parts (\ref{Poisson}), we arrive at following estimates for $|x|>r$:
\begin{equation}\label{est_w}
|w^D(x,y)|\le \frac{{C(n,m)}~\! y^{2\sigma}~\!\|u\|_{L_1(\Omega)}}
{\left((|x|-r)^2+y^2\right)^{\frac{n+m+\sigma}{2}}};
\qquad 
|\nabla w^D(x,y)|\le \frac{{C(n,m)}~\! y^{2\sigma-1}~\!\|u\|_{L_1(\Omega)}}
{\left((|x|-r)^2+y^2\right)^{\frac{n+m+\sigma}{2}}}~\!.
\end{equation}
Following \cite[Theorem 3]{FL}, we define, for $x\in \overline B_R$ and $y\ge0$, the function
$$
\widetilde w(x,y)= w^D(x,y)-\widetilde \phi(x,y),
$$
where $\widetilde \phi(\cdot,y)$ is the harmonic extension of $w^D(\cdot,y)$ on $B_R$, that is,
$$
-\Delta_x \widetilde \phi(\cdot,y)=0\quad\text{in} \ B_R;\qquad \widetilde \phi(\cdot,y)=w^D(\cdot,y)\quad\text{on} \ {\mathbb S}_R.
$$
Clearly, $\widetilde w\big|_{y=0}=(-\Delta)^ku$ and $\widetilde w\big|_{x\in {\mathbb S}_R}=0$. Further, 
we have
\begin{eqnarray}
\label{difference}
&&\int\limits_0^\infty\int\limits_{B_R}  y^{1-2\sigma}|\nabla\widetilde w|^2\,dxdy
=\int\limits_0^\infty\int\limits_{B_R}  
y^{1-2\sigma}(|\nabla w^D|^2-2\nabla w^D\cdot\nabla \widetilde \phi+|\nabla\widetilde \phi|^2)\,dxdy
\nonumber\\
&&=\int\limits_0^\infty\int\limits_{B_R} y^{1-2\sigma}|\nabla w^D|^2\,dxdy
-2{\int\limits_0^\infty\int\limits_{{\mathbb S}_R} y^{1-2\sigma}(\nabla w^D\!\cdot{\bf n})~\!
\widetilde \phi\,d{\mathbb S}_R(x)dy}
\nonumber
\\
&&\hphantom{=\int\limits_0^\infty\int\limits_{B_R} y^{1-2\sigma}|\nabla w^D|^2\,dxdy}
+{\int\limits_0^\infty\int\limits_{B_R} y^{1-2\sigma}|\nabla\widetilde \phi(x,y)|^2\,dxdy}.
\end{eqnarray}
{Since $\widetilde \phi(\cdot ,y)=w^D(\cdot,y)$ on ${\mathbb S}_R$, we can use (\ref{est_w}) to get 
$$
\Big|\int\limits_0^\infty\int\limits_{{\mathbb S}_R} y^{1-2\sigma}(\nabla w^D\!\cdot{\bf n})~\!
\widetilde \phi\,d{\mathbb S}_R(x)dy\Big|\le
\frac {{C(n,m)}~\! R^{n-1}}{(R-r)^{2n+2m-1}}\cdot\|u\|^2_{L_1(\Omega)}.
$$}
Now we estimate the last integral in (\ref{difference}). It is easy to see that
$|\nabla\widetilde \phi(\cdot,y)|^2$ is subharmonic in $B_R$ and thus the function
$$
\rho\mapsto\frac{1}{\rho^{n-1}}\int\limits_{{\mathbb S}_\rho}|\nabla\widetilde \phi(x,y)|^2\, d{\mathbb S}_\rho(x)
$$
is nondecreasing for $\rho\in(0,R)$. This implies
\begin{eqnarray*}
\int\limits_{B_R}|\nabla\widetilde \phi(x,y)|^2\,dx&=&
\int\limits_0^R\int\limits_{{\mathbb S}_{\rho}} |\nabla\widetilde \phi(x,y)|^2\,d{\mathbb S}_\rho(x)d\rho\\
&\le&\frac{R}{n}~\!\int\limits_{{\mathbb S}_R} 
(|\nabla_x\widetilde \phi(x,y)|^2+|\partial_y\widetilde \phi(x,y)|^2)\,d{\mathbb S}_R(x).
\end{eqnarray*}
Using the fact that $\partial_y\widetilde \phi(x,y)=\partial_yw^D(x,y)$ for $x\in{\mathbb S}_R$ and the well known 
estimate
$$
\int\limits_{{\mathbb S}_R} |\nabla_x\widetilde \phi(x,y)|^2\,d{\mathbb S}_R(x)
\le C(n)\int\limits_{{\mathbb S}_R} |\nabla_xw^D(x,y)|^2\,d{\mathbb S}_R(x),
$$
we can apply (\ref{est_w}) and arrive at
$$
\int\limits_0^\infty\int\limits_{B_R} y^{1-2\sigma}|\nabla\widetilde \phi(x,y)|^2\,dxdy\le
\frac {{C(n,m)}~\! R^n}{(R-r)^{2n+2m}}\cdot\|u\|^2_{L_1(\Omega)}.
$$
In conclusion, from (\ref{difference}) we infer
\begin{equation}\label{difference2}
\int\limits_0^\infty\int\limits_{B_R}  y^{1-2\sigma}|\nabla\widetilde w|^2\,dxdy
\le\int\limits_0^\infty\int\limits_{B_R} y^{1-2\sigma}|\nabla w^D|^2\,dxdy
+\frac {{C(n,m)}~\! R^n}{(R-r)^{2n+2m}}\cdot\|u\|^2_{L_1(\Omega)}.
\end{equation}

Now we use the Stinga--Torrea characterization of $Q_{\sigma,\Omega}^N$. 
Namely, a quite general result of \cite{ST} implies that
\begin{equation}\label{QN}
Q_{m,\Omega}^N[u]=Q_{\sigma,\Omega}^N[(-\Delta)^ku]=c_2(n,\sigma)\!\!
\inf\limits_{\scriptstyle w|_{x\in\partial\Omega}=0\atop \scriptstyle w|_{y=0}=(-\Delta)^ku}\,
\int\limits_0^\infty\int\limits_{\Omega} y^{1-2\sigma}|\nabla w|^2\,dxdy.
\end{equation}
Relations (\ref{QN}), (\ref{difference2}) and (\ref{QD}) give us
\begin{multline*}
Q_{m,\Omega}^N[u]\le Q_{m,B_R}^N[u]\le 
c_2(n,\sigma)~\!\int\limits_0^\infty\int\limits_{B_R} y^{1-2\sigma}|\nabla \widetilde w|^2\,dxdy
\\
\le c_2(n,\sigma)~\!\int\limits_0^\infty\int\limits_{B_R} y^{1-2\sigma}|\nabla w^D|^2\,dxdy+
\frac {{C(n,m)}~\! R^n}{(R-r)^{2n+2m}}\cdot\|u\|^2_{L_1(\Omega)}
\\
\le Q_{m,\Omega}^D[u]+\frac {{C(n,m)}~\! R^n}{(R-r)^{2n+2m}}\cdot\|u\|^2_{L_1(\Omega)},
\end{multline*}
and (\ref{new1}) follows.\smallskip

\smallskip

\noindent 
{\bf Proof of (\ref{new2})}. Let $m=2k-\sigma$, $k\in\mathbb N$, $\sigma\in(0,1)$. Denote by $w^{-D}(x,y)$, $x\in\mathbb{R}^n$, 
$y>0$, the ``dual'' Caffarelli--Silvestre extension of $(-\Delta)^ku$ (see \cite{CbSr} and \cite{FL2}), 
that is the solution of the boundary value problem
$$-\div (y^{1-2\sigma}\nabla w)=0\quad \mbox{in}\quad \mathbb R^n\times\mathbb R_+;
\qquad y^{1-2\sigma}\partial_yw\big|_{y=0}=-(-\Delta)^ku,
$$
given by the formula
\begin{equation}\label{Poisson1}
w^{-D}(x,y)=c_3(n,\sigma) \irn\frac{(-\Delta)^k u(\xi)}
{\left(|x-\xi|^2+{y^2}\right)^{\frac{n-2\sigma}{2}}}\,d\xi .
\end{equation}
Note that the representation (\ref{Poisson1}) is true also for $n=1<2\sigma$ while for $n=1$, $\sigma=1/2$
it should be rewritten as follows:
$$	
w^{-D}(x,y)=c_3(1,1/2) \irn(-\Delta)^k u(\xi)\ln(|x-\xi|^2+{y^2})\,d\xi .
$$
It is also shown in \cite{FL2} that
\begin{eqnarray}\label{Q-D}
Q_{m,\Omega}^D[u]&=&Q_{-\sigma,\Omega}^D[(-\Delta)^ku]
\\
&=&\frac 1{c_2(n,\sigma)}~\!
\bigg(2\int\limits_{\mathbb R^n} (-\Delta)^ku(x)w^{-D}(x,0)\,dx-
\int\limits_0^\infty\int\limits_{\mathbb R^n} y^{1-2\sigma}|\nabla w^{-D}|^2\,dxdy\bigg).
\nonumber
\end{eqnarray}
Integrating by parts (\ref{Poisson1}), we arrive at following estimates for $|x|>r$:
\begin{equation}\label{est_w_min}
|w^{-D}(x,y)|\le \frac{{C(n,m)}~\! \|u\|_{L_1(\Omega)}}
{\left((|x|-r)^2+y^2\right)^{\frac{n+m-\sigma}{2}}};
\qquad 
|\nabla w^{-D}(x,y)|\le \frac{{C(n,m)}~\! \|u\|_{L_1(\Omega)}}
{\left((|x|-r)^2+y^2\right)^{\frac{n+m+1-\sigma}{2}}}.
\end{equation}
Now we define,  as in  \cite[Theorem 2]{FL2}, 
$$
\widehat w(x,y)= w^{-D}(x,y)-\widehat\phi(x,y),\qquad x\in \overline B_R,\ y\ge0,
$$
where 
$$
-\Delta_x \widehat\phi(\cdot,y)=0\quad\text{in} \ B_R;
\qquad \widehat\phi(\cdot,y)=w^{-D}(\cdot,y)\quad\text{on} \ {\mathbb S}_R.
$$
Clearly, $\widehat w\big|_{x\in {\mathbb S}_R}=0$. 
Arguing as for (\ref{new1}) and using (\ref{est_w_min}) instead  of (\ref{est_w}), we obtain
\begin{equation}\label{difference1}
\int\limits_0^\infty\int\limits_{B_R}  y^{1-2\sigma}|\nabla\widehat w|^2\,dxdy
\le\int\limits_0^\infty\int\limits_{B_R} y^{1-2\sigma}|\nabla w^{-D}|^2\,dxdy
+\frac {{C(n,m)}~\! R^n}{(R-r)^{2n+2m}}\cdot\|u\|^2_{L_1(\Omega)}.
\end{equation}
We can use the ``dual'' Stinga--Torrea characterization of $Q_{-\sigma,\Omega}^N$. 
It was proved in \cite{FL2} that
\begin{eqnarray}\label{Q-N}
Q_{m,\Omega}^N[u] &=&Q_{-\sigma,\Omega}^N[(-\Delta)^ku]
\\
&= &\frac 1{c_2(n,\sigma)}
\sup\limits_{\scriptstyle w|_{x\in\partial\Omega}=0}
\bigg(\,\int\limits_{\Omega} (-\Delta)^ku(x)w(x,0)\,dx-
\int\limits_0^\infty\int\limits_{\Omega} y^{1-2\sigma}|\nabla w|^2\,dxdy\bigg).
\nonumber
\end{eqnarray}
Relations (\ref{Q-N}), (\ref{difference1}), (\ref{Q-D}) and the evident equality
$$
\int\limits_{B_R} (-\Delta)^ku(x)\widehat\phi(x,0)\,dx=0~\!,
$$
give us
\begin{multline*}
Q_{m,\Omega}^N[u]\ge Q_{m,B_R}^N[u]\ge 
\frac 1{c_2(n,\sigma)}~\!\bigg(2\int\limits_{B_R} (-\Delta)^ku(x)\widehat w(x,0)\,dx-
\int\limits_0^\infty\int\limits_{B_R} y^{1-2\sigma}|\nabla \widehat w|^2\,dxdy\bigg)
\\
\ge \frac 1{c_2(n,\sigma)}~\!\bigg(2\int\limits_{B_R} (-\Delta)^ku(x)w^{-D}(x,0)\,dx-
\int\limits_0^\infty\int\limits_{B_R} y^{1-2\sigma}|\nabla w^{-D}|^2\,dxdy\bigg)
\\
-\frac {{C(n,m)}~\! R^n}{(R-r)^{2n+2m}}\cdot\|u\|^2_{L_1(\Omega)}
\le Q_{m,\Omega}^D[u]-\frac {{C(n,m)}~\! R^n}{(R-r)^{2n+2m}}\cdot\|u\|^2_{L_1(\Omega)},
\end{multline*}
and (\ref{new2}) follows. The proof is complete.
\QED

\section{The Brezis--Nirenberg effect for Navier fractional Laplacians}\label{BNN}

We recall the Sobolev and Hardy inequalities
\begin{eqnarray}
\label{eq:Sobolev}
Q_m[u]&\ge& {\cal S}_m\,\bigg(\,\int\limits_{\R^n}|u|^\mstar\,dx\bigg)^{2/\mstar}\\
\label{eq:Hardy}
Q_m[u]&\ge& {\cal H}_m\int\limits_{\R^n}|x|^{-2m}|u|^2\,dx~\!,
\end{eqnarray}
that hold for any $u\in  {\cal C}^\infty_0(\R^n)$ and $0<m<\frac n2$. 
The {\em best Sobolev constant} ${\cal S}_m$ and the {\em best Hardy constant} ${\cal H}_m$
were explicitly computed in \cite{CoTa} and in \cite{He}, respectively.

It is well known that ${\cal H}_m$ is not attained, that is, there are no functions with finite left- 
and right-hand sides of (\ref{eq:Hardy}) providing equality in (\ref{eq:Hardy}). In contrast, it has 
been proved in \cite{CoTa} that ${\cal S}_m$ is attained by a unique family of functions, all of them
being obtained from
\begin{equation}
\label{eq:AT}
\phi(x)=(1+|x|^2)^{\frac{2m-n}{2}}
\end{equation}
by translations, dilations in $\R^n$ and multiplication by constants.

A standard dilation argument implies that 
$$
\inf_{\scriptstyle u\in{\rm Dom}(Q_{m,\Omega}^D)\atop \scriptstyle u\neq 0}\frac
{Q_{m,\Omega}^D[u]}{
\Big(\,\int\limits_{\Omega}|u|^\mstar\,dx\Big)^{2/\mstar}}
={\cal S}_m.
$$
The key fact used in further considerations is the equality
\begin{equation}
\label{eq:Sobolev_constant}
\inf_{\scriptstyle u\in{\rm Dom}(Q_{m,\Omega}^N)\atop \scriptstyle u\neq 0}\frac
{Q_{m,\Omega}^N[u]}{
\Big(\,\int\limits_{\Omega}|u|^\mstar\,dx\Big)^{2/\mstar}}
={\cal S}_m~\!,
\end{equation}
{that has been}
established in \cite{MN-Mit} (see also earlier results \cite{Ge, VdV} for $m=2$, 
\cite{GGS} for $m\in\mathbb N$ and \cite{FL} for $0<m<1$).
Clearly, the Sobolev constant ${\cal S}_m$ is never achieved on ${\rm Dom}(Q_{m,\Omega}^N)$.

The corresponding equality for the Hardy constant, that is,
\begin{equation}
\label{eq:Hardy_constant}
\inf_{\scriptstyle u\in{\rm Dom}(Q_{m,\Omega}^N)\atop \scriptstyle u\neq 0}\frac
{Q_{m,\Omega}^N[u]}
{\int\limits_{\Omega}^{\vphantom{.}}|x|^{-2m}|u|^2\,dx}={\cal H}_m~\!,
\end{equation}
was  proved in \cite{MN-Mit} as well (see also \cite{Mit2} and \cite{GGM} for $m\in\mathbb N$).

We point out that the infima
\begin{equation}
\label{eq:Poincare}
\Lambda_1(m,s):=\inf_{\scriptstyle u\in {\rm Dom}(Q_{m,\Omega}^N)\atop\scriptstyle u\ne 0}
\frac{ Q_{m,\Omega}^N[u]}{ Q_{s,\Omega}^N[u]}~~,\qquad
\widetilde\Lambda_1(m,s):=\inf_{\scriptstyle u\in {\rm Dom}(Q_{m,\Omega}^N[u])\atop\scriptstyle u\ne 0}
\frac{ Q_{m,\Omega}^N[u]}
{\int\limits_{\Omega}^{\vphantom{.}}|x|^{-2s}|u|^2\,dx}
\end{equation}
are positive and achieved. Since ${\rm Dom}(Q_{m,\Omega}^N)$ is compactly embedded into 
${\rm Dom}(Q_{s,\Omega}^N)$, this fact is well known for $\Lambda_1(m,s)$ and follows from
(\ref{eq:Hardy_constant}) for $\widetilde\Lambda_1(m,s)$. \medskip



{Weak} solutions to (\ref{eq:problem1}), (\ref{eq:problem2}) can be obtained 
as suitably normalized critical points of the functionals
\begin{gather}\label{eq:functional1}
{\cal R}^\Omega_{\lambda,m,s}[u]=
\frac
{Q_{m,\Omega}^N[u]-\lambda Q_{s,\Omega}^N[u]}
{\Big(\,\int\limits_{\Omega}|u|^\mstar\,dx\Big)^{2/\mstar}}\,,\\
~\nonumber\\
\label{eq:functional2}
\widetilde {\cal R}^\Omega_{\lambda,m,s}[u]=
\frac
{Q_{m,\Omega}^N[u]-\lambda
\int\limits_{\Omega}|x|^{-2s}|u|^2\,dx}
{\Big(\,\int\limits_{\Omega}|u|^\mstar\,dx\Big)^{2/\mstar}}\,,
\end{gather}
respectively. It is easy to see that both functionals 
are well defined on ${\rm Dom}(Q_{m,\Omega}^N)\setminus\{0\}$.

In fact, we prove the existence of ground states for functionals (\ref{eq:functional1}) and 
(\ref{eq:functional2}). We introduce the quantities
\begin{equation*}
{\cal S}^\Omega_\lambda(m,s)=
\inf_{\scriptstyle u\in {\rm Dom}(Q_{m,\Omega}^N)\atop\scriptstyle u\ne 0}{\cal R}^\Omega_{\lambda,m,s}[u];
\qquad \widetilde {\cal S}^\Omega_\lambda(m,s)=
\inf_{\scriptstyle u\in {\rm Dom}(Q_{m,\Omega}^N)\atop\scriptstyle u\ne 0}\widetilde {\cal R}^\Omega_{\lambda,m,s}[u]\,.
\end{equation*}
By standard arguments we have ${\cal S}^\Omega_\lambda(m,s) \le {\cal S}_m$. 
In addition, if $\lambda\le 0$ then ${\cal S}^\Omega_\lambda(m,s)= {\cal S}_m$ 
and it is not achieved. Similar statements hold for $\widetilde {\cal S}^\Omega_\lambda(m,s)$.

We are in position to prove our existence result that includes Theorem \ref{main1}
in the introduction.

\begin{Theorem}
\label{T:main}
Assume $s\ge 2m-\frac{n}{2}$. 

i) For any $0<\lambda<\Lambda_1(m,s)$ the infimum  ${\cal S}^\Omega_\lambda(m,s)$ is achieved and
(\ref{eq:problem1}) has a nontrivial solution in ${\rm Dom}(Q_{m,\Omega}^N)$.

ii) For any $0<\lambda<\widetilde\Lambda_1(m,s)$ the infimum $\widetilde {\cal S}^\Omega_\lambda(m,s)$ 
is achieved and (\ref{eq:problem2}) has a nontrivial solution in ${\rm Dom}(Q_{m,\Omega}^N)$.
\end{Theorem}

\proof
We prove $i)$, the proof of the second statement is similar.
Using the relation (\ref{eq:Sobolev_constant}) and arguing for instance as in \cite{MN-BN} one has that if 
$0<{\cal S}^\Omega_\lambda(m,s)< {\cal S}_m$, then ${\cal S}^\Omega_\lambda(m,s)$ is achieved. 
 
Since $0<\lambda<\Lambda_1(m,s)$,  then ${\cal S}^\Omega_\lambda(m,s)>0$ by (\ref{eq:Poincare}).

To obtain the strict inequality ${\cal S}^\Omega_\lambda(m,s)< {\cal S}_m$ we follow  \cite{BN}, and
we take advantage of the computations in \cite{MN-BN}.

Let $\phi$ be the extremal of the Sobolev inequality (\ref{eq:Sobolev}) given by 
(\ref{eq:AT}). In particular, 
\begin{equation}
\label{eq:M}
M:=Q_m[\phi]={\cal S}_m
\bigg(\,\int\limits_{\R^n}|\phi|^{\mstar}~dx\bigg)^{2/\mstar}.
\end{equation}
Fix a cutoff function $\f\in {\cal C}^\infty_0(\Omega)$,
such that $\f\equiv 1$ on the  ball $\{|x|<\delta\}$ and $\f\equiv 0$ outside the ball
$\{|x|<2\delta\}$. 

If  $\eps>0$ is small enough, the function
$$
u_{\eps}(x):=\eps^{2m-n}\f(x)\phi\left({\frac{x}{\eps}}\right)=\f(x)\left(\eps^2+|x|^2\right)^{\frac{2m-n}{2}}
$$
has compact support in $\Omega$.

From \cite[Lemma 3.1]{MN-BN} we conclude
$$
\begin{array}{lcl}
\label{eq:Hm_estimate}
\mathfrak A^\eps_m:= Q_{m,\Omega}^D[u_{\eps}]&\le&\eps^{2m-n}\left(M+C(\delta)~\!\eps^{n-2m}\right)\\
~\\
{\cal A}^\eps_s:= \displaystyle\int\limits_{\Omega}|x|^{-2s}|u_{\eps}|^2\,dx&\ge&
\begin{cases}
\displaystyle{C(\delta)~\!\eps^{4m-n-2s} }&\text{if}\quad s>2m-\frac{n}{2}\\
\\
\displaystyle{C(\delta)~\!|\log\eps|}&\text{if}\quad s=2m-\frac{n}{2}
\end{cases}
\\
~\\
\widetilde{\mathfrak A}^\eps_s:= Q_{s,\Omega}^N[u_{\eps}]
&\ge& {{\cal H}_s}~\!{\cal A}^\eps_s\qquad\qquad  \text{[~see (\ref{eq:Hardy_constant})~]}\\
~\\
{\cal B}^\eps:= \displaystyle\int\limits_{\Omega}|u_{\eps}|^\mstar\,dx&\ge&\eps^{-n}\left(({M}{\cal S}_m^{-1})^{\mstar/2}-C(\delta)~\!\eps^n\right)~\!.
\end{array}
$$
If $m$ is an integer or if $\lfloor m\rfloor\!{\not\vdots}\, 2$, then by (\ref{old2})
$$
\widetilde{\mathfrak A}^\eps_m:= Q_{m,\Omega}^N[u_{\eps}]\le \mathfrak A^\eps_m,
$$
and we obtain
\begin{eqnarray}\label{eq:ups1}
&&{\cal R}^\Omega_{\lambda,m,s}[u_\eps]\le
{\cal S}_m~\!\frac {1+C(\delta)~\!\eps^{n-2m}-\lambda C(\delta)~\!\eps^{2m-2s}}
{1-C(\delta)~\!\eps^n}\,, \qquad\text{if}\quad s>2m-\frac{n}{2}\\
\nonumber\\
\label{eq:ups2}
&&{\cal R}^\Omega_{\lambda,m,s}[u_\eps]\le
{\cal S}_m~\!\frac {1+C(\delta)~\!\eps^{n-2m}-\lambda C(\delta)~\!\eps^{n-2m}|\log\eps|}
{1-C(\delta)~\!\eps^n}\,, \quad\text{if}\ \ s=2m-\frac{n}{2}.
\end{eqnarray}
Thus, for any sufficiently small $\eps>0$ we have that ${\cal R}^\Omega_{\lambda,m,s}[u_\eps]< {\cal S}_m$, and
the statement follows.\medskip

It remains to consider the case $\lfloor m\rfloor\vdots\, 2$. Since $\|u_\eps\|_{L_1(\Omega)}\le C(\delta)$,
the estimate (\ref{new1}) implies
$$
\widetilde{\mathfrak A}^\eps_m\le \mathfrak A^\eps_m+C(\delta)=
\eps^{2m-n}\left(M+C(\delta)~\!\eps^{n-2m}\right),
$$
and we again arrive at (\ref{eq:ups1}), (\ref{eq:ups2}).
\QED

For the case $s< 2m-\frac{n}{2}$ we limit ourselves to point out
the next simple existence result, as in \cite{MN-BN}.

\begin{Theorem}
\label{T:critical}
Assume $s< 2m-\frac{n}{2}$. 
\begin{description}
\item$~i)$ There exists $\lambda^*\in[0,\Lambda_1(m,s))$ such that
for any $\lambda\in(\lambda^*,\Lambda_1(m,s))$ the infimum ${\cal S}^\Omega_\lambda(m,s)$ is attained, and
hence (\ref{eq:problem1}) has a nontrivial solution.
\item$ii)$ There exists $\widetilde\lambda^*\in[0,\widetilde\Lambda_1(m,s))$ such that for any 
$\lambda\in(\widetilde\lambda^*,\widetilde\Lambda_1(m,s))$ the infimum $\widetilde{\cal S}^\Omega_\lambda(m,s)$ 
is attained, and hence (\ref{eq:problem2}) has a nontrivial solution.
\end{description}
\end{Theorem}

\footnotesize
\label{References}

\end{document}